\documentclass[12pt,oneside,reqno]{amsart}
\usepackage{graphicx}
\usepackage{mathrsfs}
\usepackage{stmaryrd}
\usepackage{amsfonts}
\usepackage{enumerate,amsmath,amssymb,amsthm}
\newcommand{\dif}{\mathrm{d}}

\newcommand{\be}{\begin{eqnarray}}
\newcommand{\ee}{\end{eqnarray}}
\newcommand{\ce}{\begin{eqnarray*}}
\newcommand{\de}{\end{eqnarray*}}
\newtheorem{theorem}{Theorem}[section]
\newtheorem{lemma}[theorem]{Lemma}
\newtheorem{remark}[theorem]{Remark}
\newtheorem{definition}[theorem]{Definition}
\newtheorem{proposition}[theorem]{Proposition}
\newtheorem{Example}[theorem]{Example}
\newtheorem{corollary}[theorem]{Corollary}
\def\e{\varepsilon}
\def\t{\theta}

\def\[{{\Big[}}
\def\]{{\Big]}}
\def\<{{\langle}}
\def\>{{\rangle}}
\def\({{\Big(}}
\def\){{\Big)}}

\def\no{\nonumber}
\def\bt{\begin{theorem}}
\def\et{\end{theorem}}
\def\bl{\begin{lemma}}
\def\el{\end{lemma}}
\def\br{\begin{remark}}
\def\er{\end{remark}}
\def\bx{\begin{Example}}
\def\ex{\end{Example}}
\def\bd{\begin{definition}}
\def\ed{\end{definition}}
\def\bp{\begin{proposition}}
\def\ep{\end{proposition}}
\def\bc{\begin{corollary}}
\def\ec{\end{corollary}}

\def\cB{{\mathcal B}}

\def\cL{{\mathcal L}}
\def\cM{{\mathcal M}}

\def\cP{{\mathcal P}}

\def\mE{{\mathbb E}}

\def\mN{{\mathbb N}}

\def\mP{{\mathbb P}}

\def\mR{{\mathbb R}}
\def\mS{{\mathbb S}}

\def\mU{{\mathbb U}}

\def\sA{{\mathscr A}}
\def\sB{{\mathscr B}}

\def\sF{{\mathscr F}}

\def\sL{{\mathscr L}}

\def\geq{\geqslant}
\def\leq{\leqslant}

\begin{document}

\allowdisplaybreaks

\title{Limit theorems of stochastic differential equations with jumps*}

\author{Huijie Qiao}

\thanks{{\it AMS Subject Classification(2020):} 60H10, 60J76}

\thanks{{\it Keywords:} Weak solutions, martingale solutions, Fokker-Planck equations, superposition principles}

\thanks{*This work was partly supported by NSF of China (No. 11001051, 11371352, 12071071) and China Scholarship Council under Grant No. 201906095034.}

\subjclass{}

\date{}

\dedicatory{School of Mathematics,
Southeast University\\
Nanjing, Jiangsu 211189,  China\\
Department of Mathematics, University of Illinois at
Urbana-Champaign\\
Urbana, IL 61801, USA\\
hjqiaogean@seu.edu.cn}

\begin{abstract} 
In this paper, we study the convergence for solutions to a sequence of (possibly degenerate) 
stochastic differential equations with jumps, when the coefficients converge in some appropriate sense. Our main tools are the superposition principles. And then we analyze some special cases and give some concrete and verifiable conditions.
\end{abstract}

\maketitle \rm

\section{Introduction}

%
Stochastic differential equations (SDEs in short) with jumps appear naturally in various applied fields. 
And more and more people pay attention to them. For example, in \cite{jj}, Jacod systematically discussed the martingale problems of SDEs with jumps. In \cite{jjas}, Jacod and Shiryaev studied limit theorems of SDEs driven by c\`adl\`ag processes under Lipschitz conditions.
Qiao and Zhang \cite{q4} proved that, under non-Lipschitz conditions, for almost
all sample points $\omega$, the solutions to a certain SDE with jumps form a homeomorphism flow.
Recently, Fournier and Xu \cite{FX} established the equivalence between SDEs with jumps and 
the corresponding non-local Fokker-Planck equations under only linear growth conditions.


In this paper, we study the convergence for solutions to a sequence of (possibly degenerate) 
stochastic differential equations with jumps, when the coefficients converge in some appropriate sense. More precisely, we fix a $T>0$ and consider the following sequence of SDEs with jumps:
\be
\dif X_t^n=b^n(t, X_t^n)\dif t+\sigma^n(t, X_t^n)\dif B_t+\gamma^n\int_{\mU}g(t,X_{t-}^n,u) N(\dif t, \dif u), \quad t\in[0,T],
\label{Eq0}
\ee
where $(B_t)$ is an $m$-dimensional Brownian motion and $N(\dif t, \dif u)$ is 
a Poisson random measure with the intensity $\dif t\nu(\dif u)$. Here $\nu$ is a finite measure defined on $(\mU,\mathscr{U})$, where $(\mU,\mathscr{U})$ is any measurable space. The coefficients $b^n: [0,T]\times\mR^d\mapsto\mR^d$, $\sigma^n: [0,T]\times\mR^d\mapsto\mR^{d\times m}$ are Borel measurable functions, $\{\gamma^n\}$ is a sequence of real numbers and $g: [0,T]\times\mR^d\times\mU\mapsto\mR^d$ is Borel measurable.  Under some pretty weak conditions, we show that, when $b^n\rightarrow b, \sigma^n\rightarrow\sigma$, $\gamma^n\rightarrow\gamma$ in some sense as $n\rightarrow\infty$, the martingale solutions of Eq.(\ref{Eq0}) converge to that of 
the following equation
\be
\dif X_t=b(t,X_t)\dif t+\sigma(t,X_t)\dif B_t+\gamma\int_{\mU} g(t, X_{t-}, u) N(\dif t, \dif u), \qquad t\in[0,T],
\label{Eq01}
\ee
where $b: [0,T]\times\mR^d\mapsto\mR^d$, $\sigma: [0,T]\times\mR^d\mapsto \mR^{d\times m}$ are Borel measurable and $\gamma$ is a real number. This kind of convergence results are very useful in approximation theory and statistics (\cite{jjas}). 
And then we analyze some special cases and give some concrete and verifiable conditions.

It is worthwhile to mentioning our conditions and methods. Here we require that the jump measure $\nu$ is finite so that some concrete and verifiable conditions are given. This is convenient in application. Besides, we state a superposition principle for the following SDE:
\be
\dif X_t=b(t,X_t)\dif t+\sigma(t,X_t)\dif B_t+\int_{\mU}f(t,X_{t-},u)N(\dif t, \dif u), \quad t\in[0,T],
\label{Eq1}
\ee
where $f: [0,T]\times\mR^d\times\mU\mapsto\mR^d$ is Borel measurable. And then we apply  the superposition principle to prove our convergence theorems. So we effectively avoid some properties, such as the ergodicity and regularity, which usually are used in the proofs of limit theorems.

In this paper, our motivation lies in offering some certain weak conditions of convergence for the martingale solutions to Eq.(\ref{Eq0}). Especially, we give out concrete and verifiable conditions under which the martingale solutions to Eq.(\ref{Eq0}) converge to that of a SDE with pure jumps (Corollary \ref{nocom}), that of a SDE without jumps (Proposition \ref{nojump}) and that of an ordinary differential equation (Proposition \ref{ODE}), respectively. Thus, it is convenient to apply these convergence theorems to  approximation theory and statistics (\cite{jjas}). 


The paper is arranged as follows. In the next section, we introduce some concepts, such as weak solutions and martingale solutions of SDEs with jumps, and weak solutions of Fokker-Planck equations, and their relationship. We study limits of SDEs with jumps in Section \ref{add}. In Section \ref{spe}, we analyze some special cases and give some concrete and verifiable conditions. Finally, we prove Remark \ref{bnan} in the appendix.

The following convention will be used throughout the paper: $C$ with
or without indices will denote different positive constants
whose values may change from one place to
another.
 
\section{Preliminary}\label{pre}

\subsection{Notation}

In this subsection, we introduce some  notation used in the sequel. 

We use $\mid\cdot\mid$ and $\parallel\cdot\parallel$  for the norms of vectors and matrices, respectively. Let $\langle\cdot$ , $\cdot\rangle$ be the scalar product in $\mR^d$. Let $A^*$ denote the transpose of the matrix $A$.

$C^2(\mR^d)$ stands for the space of continuous functions on $\mR^d$ which have continuous partial derivatives of order up to $2$, and $C_b^2(\mR^d)$ stands for the subspace of $C^2(\mR^d)$, consisting of functions whose derivatives up to order 2 are bounded.
$C_c^2(\mR^d)$ is the collection of all functions in $C^2(\mR^d)$ with compact supports and $C_c^\infty(\mR^d)$ denotes the collection of all real-valued $C^\infty$ functions of compact
supports.

Let $\sB(\mR^d)$ be the Borel $\sigma$-field on $\mR^d$. Let $\cP({\mR^d})$ be the space of all probability measures on $\sB(\mR^d)$, equipped with the topology of weak convergence. 
Let $\cP_1(\mR^d)$ be the collection of all the probability measures $\mu$ on $\sB(\mR^d)$ satisfying
\ce
\mu(|\cdot|):=\int_{\mR^d}\mid{x}\mid\,\mu(\dif x)<\infty.
\de
Let $L^\infty([0,T],\cP_1(\mR^d))$ be the collection of all measurable families $(\mu_t)_{t\in[0,T]}$ of probability measures on $\sB(\mR^d)$ satisfying $\sup\limits_{t\in[0,T]}\mu_t(|\cdot|)<\infty$.

\subsection{Weak solutions and martingale solutions for SDEs with jumps}
In this subsection, we introduce the concepts of weak solutions and martingale solutions for SDEs with jumps, and study their relationship.

First of all, we recall the definition of weak solutions to Eq.(\ref{Eq1}). \cite{ks} is a good reference for the definition below, and it does not deal with SDEs with jumps.

\bd(Weak solutions)\label{weaksolu}
By a weak solution  to Eq.(\ref{Eq1}), we mean a septuple $\{(\Omega,\sF,\mP;(\sF_t)_{t\in[0,T]}), (B,N,X)\}$, where $(\Omega,\sF,\mP;(\sF_t)_{t\in[0,T]})$ is a complete filtered probability space, $(B_t)$ is an $(\sF_t)$-adapted Brownian motion, $N(\dif t, \dif u)$ is an 
$(\sF_t)$-adapted Poisson random measure, independent of $(B_t)$, with the intensity $\dif t\nu(\dif u)$, and $(X_t)$ is an $(\sF_t)$-adapted process such that  for all $t\in[0,T]$,
\ce
\mP\left(\int_0^t\left(|b(s,X_s)|+\|\sigma\sigma^*(s,X_s)\|+\int_{\mU}|f(s,X_{s-},u)|\nu(\dif u)\right)\dif s<\infty\right)=1,
\de 
and
\ce
X_t=X_0+\int_0^tb(s,X_s)\dif s+\int_0^t\sigma(s,X_s)\dif B_s+\int_0^t\int_{\mU}f(s,X_{s-},u)N(\dif s, \dif u), a.s.\mP.
\de

If two weak solutions to Eq.(\ref{Eq1}), $\{(\Omega,\sF,\mP;(\sF_t)_{t\in[0,T]}), (B,N,X^1)\}$ and $\{(\Omega,\sF,\mP;(\sF_t)_{t\in[0,T]}), \\(B,N,X^2)\}$  with $X^1_0=X^2_0$, satisfy
$$
X^1_t=X^2_t, \quad t\in[0,T], ~ a.s.\mP,
$$
then we say pathwise uniqueness holds for Eq.(\ref{Eq1}).

If any two weak solutions to Eq.(\ref{Eq1}) with the same initial distribution have the same law, then we say uniqueness in law holds for Eq.(\ref{Eq1}).
\ed

It is known that the pathwise uniqueness implies the uniqueness in law for Eq.(\ref{Eq1}).

Let $D_T:=D([0,T], \mR^d)$ be the collection of c\`adl\`ag functions from $[0, T]$ to $\mR^d$.
The generic element in $D_T$ is denoted by $w$. 
We equip $D_T$ with the Skorokhod topology and then $D_T$ is a Polish space. For any $t\in[0,T]$, set 
$$
e_t: D_T\mapsto\mR^d, \quad e_t(w)=w_t, \quad w\in D_T.
$$
Let $\cB_t:=\sigma\{w_s: s\in[0,t]\}$,  $\bar{\cB}_t:=\cap_{s>t}\cB_s$, and $\cB:=\cB_T$. 
For $\phi\in C_b^2(\mR^d)$, set
\ce
&&(\sA_t\phi)(x):=b_i(t,x)\partial_i\phi(x)+a_{ij}(t,x)\partial_{ij}\phi(x),\no\\
&&(\sB_t\phi)(x):=\int_{\mU}\[\phi(x+f(t,x,u))-\phi(x)\]\nu(\dif u),
\de
where $a(t,x)=\frac{1}{2}\sigma\sigma^*(t,x)$. In the following, we define martingale solutions of Eq.(\ref{Eq1}).(c.f.\cite{ks, sv})

\bd(Martingale solutions)\label{martsolu}
For $\mu_0\in\cP(\mR^d)$. A probability measure $\mP$ on $(D_T, \cB)$ is called a martingale solution of Eq.(\ref{Eq1}) with the initial law $\mu_0$ at time $0$, if

(i) $\mP\circ e^{-1}_0=\mu_0$,

(ii) For any $\phi\in{C_c^2(\mR^d)}$,
\be
\cM_t^\phi&:=&\phi(w_t)-\phi(w_0)-\int_0^t(\sA_s\phi+\sB_s\phi)(w_s)ds
\label{eq2}
\ee
is a $(\bar{\cB}_t)_{t\in[0,T]}$-adapted martingale under the probability measure $\mP$.
And the uniqueness of the martingale solutions to Eq.(\ref{Eq1}) means that, for any $s\in[0,T]$ and any $\mu_s\in\cP(\mR^d)$, if $\hat{\mP}, \tilde{\mP}$ are two martingale solutions to Eq.(\ref{Eq1}) with the initial law $\mu_s$ at time $s$, then $\hat{\mP}=\tilde{\mP}$.
\ed

Next, we assume:
\begin{enumerate}[(${\bf H}_{b,\sigma}$)]
\item 
There is a constant $C_1$ such that for all $(t,x)\in[0,T]\times\mR^d$,
\ce
|b(t,x)|+\|\sigma(t,x)\|\leq C_1(1+|x|).
\de
\end{enumerate}
\begin{enumerate}[(${\bf H}_{f}$)]
\item 
There is a constant $C_2$ such that for all $(t,x)\in[0,T]\times\mR^d$,
\ce
\int_{\mU}|f(t,x,u)|^2\nu(\dif u)\leq C_2(1+|x|)^2.
\de
\end{enumerate}

By the H\"older inequality, it is easy to see that 
\begin{enumerate}[(${\bf H}^{\prime}_{f}$)]
\item 
$$
\int_{\mU}|f(t,x,u)|\nu(\dif u)\leq C(1+|x|).
$$
\end{enumerate}
The relationship between martingale solutions and weak solutions is as follows. 

\bp\label{martweak}
Assume that (${\bf H}_{b,\sigma}$) and (${\bf H}_f$) hold.

(i) For any $\mu_0\in\cP(\mR^d)$, the existence of a weak solution $(X_t)_{t\in[0,T]}$ to Eq.(\ref{Eq1}) with $\cL_{X_0}=\mu_0$ is equivalent to the existence of a martingale solution $\mP$ to Eq.(\ref{Eq1}) with the initial law $\mu_0$. Moreover, $\cL_{X_t}=\mP\circ e_t^{-1}$ for any $t\in[0,T]$.

(ii) The uniqueness of martingale solutions $\mP$ to Eq.(\ref{Eq1}) implies the uniqueness in law of weak solutions $(X_t)_{t\in[0,T]}$ to Eq.(\ref{Eq1}). 
\ep
\begin{proof}
We only prove (i). Assume that $\{(\hat{\Omega},\hat{\sF},\hat{\mP};(\hat{\sF}_t)_{t\in[0,T]}), (\hat{B},\hat{N},\hat{X})\}$ is a weak solution of Eq.(\ref{Eq1}) with $\cL_{\hat{X}_0}=\mu_0$. It follows from the It\^o formula, that for any $\phi\in{C_c^2(\mR^d)}$,
\ce
&&\phi(\hat{X}_t)-\phi(\hat{X}_0)-\int_0^t(\sA_s\phi+\sB_s\phi)(\hat{X}_s)ds\\
&=&\int_0^t\partial_i\phi(\hat{X}_s)\sigma_{ij}(s,\hat{X}_s)\dif \hat{B}_s+\int_0^t\int_{\mU}(\phi(\hat{X}_{s-}+f(s,\hat{X}_{s-},u))-\phi(\hat{X}_{s-}))\tilde{\hat{N}}(\dif s\dif u),
\de
where $\tilde{\hat{N}}(\dif s\dif u):=\hat{N}(\dif s\dif u)-\nu(\dif u)\dif s$ is the compensated martingale measure of $\hat{N}(\dif s\dif u)$. Note that 
\ce
\int_0^T|\partial_i\phi(\hat{X}_s)\sigma_{ij}(s,\hat{X}_s)|^2\dif s\leq\int_0^TCI_{|\hat{X}_s|\leq M}(1+|\hat{X}_s|)^2\dif s\leq CT(1+M)^2
\de
and
\ce
&&\int_0^T\int_{\mU}|\phi(\hat{X}_s+f(s,\hat{X}_s,u))-\phi(\hat{X}_s)|^2\nu(\dif u)\dif s\\
&\leq&\|\phi\|^2_{C_c^2(\mR^d)}\int_0^T\int_{\mU}I_{|\hat{X}_s|\leq M}|f(s,\hat{X}_s,u)|^2\nu(\dif u)\dif s\\
&&+\int_0^T\int_{\mU}I_{|\hat{X}_s|>M}|\phi(\hat{X}_s+f(s,\hat{X}_s,u))|^2\nu(\dif u)\dif s\\
&\leq&\|\phi\|^2_{C_c^2(\mR^d)}\int_0^TI_{|\hat{X}_s|\leq M}(1+|\hat{X}_s|)^2\dif s\\
&&+\|\phi\|^2_{C_c^2(\mR^d)}\nu(\mU)T\\
&\leq&\|\phi\|^2_{C_c^2(\mR^d)}T(1+M)^2+\|\phi\|^2_{C_c^2(\mR^d)}\nu(\mU)T,
\de
where $M>0$ is a number such that ${\rm supp} (\phi)\subset B_M:=\{y\in\mR^d; |y|\leq M\}$.
Thus,
$$
\int_0^t\partial_i\phi(\hat{X}_s)\sigma_{ij}(s,\hat{X}_s)\dif \hat{B}_s+\int_0^t\int_{\mU}(\phi(\hat{X}_{s-}+f(s,\hat{X}_{s-},u))-\phi(\hat{X}_{s-}))\tilde{\hat{N}}(\dif s\dif u)
$$ 
is an $(\hat{\sF}_t)_{t\in[0,T]}$-adapted martingale and then 
$$
\phi(\hat{X}_t)-\phi(\hat{X}_0)-\int_0^t(\sA_s\phi+\sB_s\phi)(\hat{X}_s)ds
$$ 
is an $(\hat{\sF}_t)_{t\in[0,T]}$-adapted martingale. 
Set $\mP:=\hat{\mP}\circ \hat{X}^{-1}_{\cdot}$, the argument above shows that $\mP$ is a martingale solution of Eq.(\ref{Eq1}). 

Conversely, assume that $\mP$ is a martingale solution of Eq.(\ref{Eq1}). 
For any $n\geq 1$, we take $\phi_n\in C_c^2(\mR^d)$ so that $\phi_n(x)=x^j$, $j$-th component
of $x$, for all $|x|\leq n$, and define $\tau_n:=\inf\{0\leq t\leq T, |w_t|\geq n\}$. It follows from (\ref{eq2}) that
\ce
\cM_t^{\phi_n}&=&\phi_n(w_{t\land\tau_n})-\phi_n(w_0)-\int_0^{t\land\tau_n}(\sA_s\phi_n+\sB_s\phi_n)(w_s)\dif s\\
&=&w^j_{t\land\tau_n}-w^j_0-\int_0^{t\land\tau_n}\(b_j(s,w_s)+\int_{\mR^d}y_j\nu_{s,w_s}(\dif y)\)\dif s
\de
is a $(\bar{\cB}_t)_{t\in[0,T]}$-adapted martingale under the probability measure $\mP$, where $\nu_{s,w_s}(\dif y):=\nu(\dif f^{-1}(s,w_s,\cdot)(y))$. Thus, $w_{t\land\tau_n}$ is a semimartingale and $w_t$ has local characteristics $(\tilde{b}, \tilde{a}, \tilde{\nu}_{\cdot,w_\cdot})$ (c.f. \cite{jj, jjas}). 
By the definition of the characteristics, we have that for any $\xi\in\mR^d$,
\be
e^{i\<\xi,w_t\>}-e^{i\<\xi,w_0\>}-\int_0^te^{i\<\xi,w_s\>}\left\{i\tilde{b}_j\xi_j-\tilde{a}_{ij}\xi_i\xi_j+\int_{\mR^d}\left(e^{i\<\xi,y\>}-1-i\<\xi,y\>\right)\tilde{\nu}_{s,w_s}(\dif y)\right\}\dif s
\label{mart1}
\ee
is a $(\bar{\cB}_t)_{t\in[0,T]}$-adapted local martingale. On the other side, applying $\psi_n(x)=e^{i\<\xi,x\>}\chi_n(x)$  to (\ref{eq2}), where $\chi_n$ is a smooth function such that $\chi_n(x)=1, |x|\leq n$ and $\chi_n(x)=0, |x|\geq 2n$, we know that 
\ce
\cM^{\psi_n}_t=\psi_n(w_t)-\psi_n(w_0)-\int_0^t(\sA_s\psi_n+\sB_s\psi_n)(w_s)\dif s
\de
is a $(\bar{\cB}_t)_{t\in[0,T]}$-adapted martingale under the probability measure $\mP$. Set $\tau_v:=\inf\{0\leq t\leq T, |w_t|>v\}$ for $v\in\mN$, and then $\{\tau_v\}$ is a $(\bar{\cB}_t)_{t\in[0,T]}$-stopping time sequence and $\tau_v\uparrow T$ as $v\rightarrow\infty$. Thus, 
\ce
\cM^{\psi_n}_{t\land\tau_v}&=&\psi_n(w_{t\land\tau_v})-\psi_n(w_0)-\int_0^{t\land\tau_v}(\sA_s\psi_n+\sB_s\psi_n)(w_s)\dif s
\de
is still a $(\bar{\cB}_t)_{t\in[0,T]}$-adapted martingale under $\mP$. By the dominated convergence theorem we obtain
\ce
\cM^{e^{i\<\xi,\cdot\>}}_{t\land\tau_v}&=&e^{i\<\xi,w_{t\land\tau_v}\>}-e^{i\<\xi,w_0\>}-\int_0^{t\land\tau_v}(\sA_se^{i\<\xi,\cdot\>}+\sB_se^{i\<\xi,\cdot\>})(w_s)\dif s
\de
is also a $(\bar{\cB}_t)_{t\in[0,T]}$-adapted martingale under $\mP$. That is,
\be
&&e^{i\<\xi,w_t\>}-e^{i\<\xi,w_0\>}-\int_0^te^{i\<\xi,w_s\>}\Big\{ib_j\xi_j+i\int_{\mR^d}y_j\xi_j\nu_{s,w_s}(\dif y)-a_{lj}\xi_l\xi_j\no\\
&&\qquad\qquad +\int_{\mR^d}(e^{i\<\xi,y\>}-1-i\<\xi,y\>)\nu_{s,w_s}(\dif y)\Big\}\dif s
\label{mart2}
\ee
is a local martingale. Comparing (\ref{mart1}) with (\ref{mart2}), one gets that $\tilde{b}(s,w_s)=b(s,w_s)+\int_{\mR^d}y\nu_{s,w_s}(\dif y)$, $\tilde{a}(s,w_s)=a(s,w_s)$ and $\tilde{\nu}_{s,w_s}=\nu_{s,w_s}$. Put
\ce
\kappa_t:=w_t-w_{t-}, \qquad N_{\kappa}((0,t],A):=\sum_{0<s\leq t}I_{A}(\kappa_s), \quad A\in\sB(\mR^d\setminus\{0\}),
\de
and then $\tilde{N}_{\kappa}((0,t],A):=N_{\kappa}((0,t],A)-\int_0^t\nu_{s,w_s}(A)\dif s$ is the compensated martingale measure of $N_{\kappa}$. Define
$$
M_t:=w_t-w_0-\int_0^t\(b(s,w_s)+\int_{\mR^d}y\nu_{s,w_s}(\dif y)\)\dif s-\int_0^t\int_{\mR^d}y\tilde{N}_{\kappa}(\dif s\dif y),
$$
and then by the definition of the characteristics, it holds that $M_t$ is a continuous local martingale with $\<M_l, M_j\>_t=\int_0^t2a_{lj}(s,w_s)\dif s$. 

Next, we take another filtered probability space $(\tilde{\Omega},\tilde{\sF},\tilde{\mP};(\tilde{\sF}_t)_{t\in[0,T]})$ and a $d$-dimensional Brownian motion $\tilde{B}$ on it. Define
\ce
(\check{\Omega},\check{\sF},\check{\mP};(\check{\sF}_t)_{t\in[0,T]}):=(\tilde{\Omega},\tilde{\sF},\tilde{\mP};(\tilde{\sF}_t)_{t\in[0,T]})\times(D_T,\cB,\mP;(\bar{\cB}_t)_{t\in[0,T]}),
\de
and 
$$
\pi:\check{\Omega}\mapsto D_T, \quad \pi(\check{\omega})=w, \quad \check{\omega}=(\tilde{\omega}, w)\in\check{\Omega}.
$$
So, on $(\check{\Omega},\check{\sF},\check{\mP};(\check{\sF}_t)_{t\in[0,T]})$, $\check{M}_t(\check{\omega}):=M_t(w)$ is still a local martingale and $\check{B}_t(\check{\omega}):=\tilde{B}_t(\tilde{\omega})$ is still a Brownian motion. By the similar method to that in \cite[4.6 Proposition, Page 315]{ks}, we can construct a $d$-dimensional Brownian motion $B$ on $(\check{\Omega},\check{\sF},\check{\mP};(\check{\sF}_t)_{t\in[0,T]})$ so that 
\ce
\check{M}_t=\int_0^t\sigma(s,\pi(\check{\omega})_s)\dif B_s,
\de
i.e.
\ce
\pi(\check{\omega})_t=\pi(\check{\omega})_0+\int_0^tb(s,\pi(\check{\omega})_s)\dif s+\int_0^t\sigma(s,\pi(\check{\omega})_s)\dif B_s+\int_0^t\int_{\mU}f(s,\pi(\check{\omega})_{s-},u)N(\dif s\dif u),
\de
where $N(\dif s\dif u)$ is a Poission random measure on $(\check{\Omega},\check{\sF},\check{\mP};(\check{\sF}_t)_{t\in[0,T]})$ with the intensity measure $\nu(\dif u)\dif s$. Therefore, $\{(\check{\Omega},\check{\sF},\check{\mP};(\check{\sF}_t)_{t\in[0,T]}),(B,N,\pi)\}$ is a weak solution of Eq.(\ref{Eq1}). Thus, the proof is complete.
\end{proof}

\subsection{Weak solutions of Fokker-Planck equations}
In this subsection, we introduce weak solutions of the Fokker-Planck equations (FPEs in short) and prove a property about them.

Consider the FPE associated with Eq.(\ref{Eq1}):
\be
\partial_t\mu_t=(\sA_t+\sB_t)^*\mu_t,
\label{FPE1}
\ee
where $(\sA_t+\sB_t)^*$ is the adjoint operator of $\sA_t+\sB_t$, and $(\mu_t)_{t\in[0,T]}$ is a family of probability measures on $\mR^d$.  Weak solutions of Eq.(\ref{FPE1}) are defined as follows. 

\bd\label{weakfpe}
A measurable family $(\mu_t)_{t\in[0,T]}$ of probability measures is called a
 weak solution of the non-local FPE (\ref{FPE1}) starting from $\mu_0$ at time $0$ if for any $R>0$ and $t\in[0,T]$,
\be
\int_0^t\int_{\mR^d}I_{B_R}(x)\left(|b(s,x)|+\|a(s,x)\|+\int_{\mU}|f(s,x,u)|\nu(\dif u)\right)\mu_s(\dif x)\dif s<\infty,
\label{deficond1}
\ee
and for all $\phi\in C_c^2(\mR^d)$ and $t\in[0,T]$,
\be
\mu_t(\phi)=\mu_0(\phi)+\int_0^t\mu_s(\sA_s\phi+\sB_s\phi)\dif s,
\label{deficond2}
\ee
where $\mu_t(\phi):=\int_{\mR^d}\phi(x)\mu_t(\dif x)$. The uniqueness of weak solutions to Eq.(\ref{FPE1}) means that, if for any $s\in[0,T]$ and any $\mu_s\in\cP(\mR^d)$, $(\hat{\mu}_t)_{t\in[s,T]}$ and $(\tilde{\mu}_t)_{t\in[s,T]}$ are two weak solutions to Eq.(\ref{FPE1}) starting from $\mu_s$ at time $s$, then $\hat{\mu}_t=\tilde{\mu}_t$ for any $t\in[s,T]$.
\ed

We claim that the definition above makes sense. That is, under (\ref{deficond1}), it holds that $\int_0^t|\mu_s(\sA_s\phi+\sB_s\phi)|\dif s<\infty$ for all $\phi\in C_c^2(\mR^d)$. Indeed, for any $\phi\in C_c^2(\mR^d)$, we assume that the support of $\phi$ is in some ball $B_R$.  Then
\ce
\int_0^t|\mu_s(\sA_s\phi+\sB_s\phi)|\dif s&\leq&\int_0^t\int_{\mR^d}\bigg[|b_i(s,x)\partial_i\phi(x)|+\frac{1}{2}|(\sigma\sigma^*)_{ij}(s,x)\partial_{ij}\phi(x)|\\
&&\qquad\qquad +\int_{\mU}|\phi(x+f(s,x,u))-\phi(x)|\nu(\dif u)\bigg]\mu_s(\dif x)\dif s\\
&\leq&\|\phi\|_{C_c^2(\mR^d)}\int_0^t\int_{\mR^d}I_{B_R}(x)\(|b(s,x)|+\|a(s,x)\|\\
&&\qquad\qquad\qquad\qquad+\int_{\mU}|f(s,x,u)|\nu(\dif u)\)\mu_s(\dif x)\dif s\\
&&+\|\phi\|_{C_c^2(\mR^d)}\int_0^t\int_{\mR^d}I_{B^c_R}(x)\nu(\mU)\mu_s(\dif x)\dif s\\
&<&\infty.
\de

If a weak solution $(\mu_t)_{t\in[0,T]}$ of the non-local FPE (\ref{FPE1}) is absolutely continuous with respect to the Lebesgue measure, then there exists a measurable family of non-negative functions $(v_t)_{t\in[0,T]}$ with $\int_{\mR^d}v_t(x)\dif x=1$ such that $\mu_t(\dif x)=v_t(x)\dif x$ for any $t\in[0,T]$. Thus, $(v_t)_{t\in[0,T]}$ satisfies the following equation in the distributional sense
\be
\partial_t v_t=-\partial_i(b_iv_t)+\partial_{ij}(a_{ij}v_t)+\sB_t^*v_t.
\label{FPE2}
\ee
Set
\ce
\sL_+&:=&\bigg\{v=(v_t)_{t\in[0,T]} : v_t\geq 0, \int_{\mR^d}v_t(x)\dif x=1, ~\mbox{for}~ \forall t\in[0,T], \\
&&\qquad\qquad ~\mbox{and}~ \sup\limits_{t\in[0,T]}\left(\int_{\mR^d}|x|v_t(x)\dif x\right)<\infty\bigg\}.
\de
If there exists a $v\in\sL_+$ satisfying Eq.(\ref{FPE2}) in the distributional sense, we say  Eq.(\ref{FPE2}) has a weak solution in $\sL_+$. Later, we will assume that Eq.(\ref{FPE2}) has a unique weak solution in $\sL_+$, that is, for any $s\in[0,T]$ 
and any non-negative measurable function $v_s$ with $\int_{\mR^d}v_s(x)\dif x=1$ and $\int_{\mR^d}|x|v_s(x)\dif x<\infty$, if $\hat{v}$ and $\tilde{v}$ are two weak solutions of Eq.(\ref{FPE2}) in $\sL_+$ starting from $v_s$ at time $s$, then $\hat{v}_t(x)=\tilde{v}_t(x), \forall x\in\mR^d$ for any $t\in[s,T]$.

\bp\label{soluinit}
Assume that (${\bf H}_{b,\sigma}$) and (${\bf H}^{\prime}_f$) hold,  and that any $\mu_0(\dif x)=\bar{v}_0(x)\dif x\in\cP_1(\mR^d)$. If there exists a measurable family of non-negative functions $(v_t)_{t\in[0,T]}$ with $\int_{\mR^d}v_t(x)\dif x=1$ and $v_0(x)=\bar{v}_0(x)$ such that $v$ is a weak solution of the non-local FPE (\ref{FPE2}), then $v\in\sL_+$.
\ep
\begin{proof}
By the assumption and the definition of $\sL_+$, we only need to prove
$$
 \sup\limits_{t\in[0,T]}\int_{\mR^d}|x|v_t(x)\dif x<\infty.
$$
By Definition \ref{weakfpe}, it holds that for all $\phi\in C_c^2(\mR^d)$ and $t\in[0,T]$
\be
\int_{\mR^d}v_t(x)\phi(x)\dif x&=&\int_{\mR^d}\bar{v}_0(x)\phi(x)\dif x+\int_0^t\int_{\mR^d}(\sA_s\phi(x)+\sB_s\phi(x))v_s(x)\dif x\dif s\no\\
&=&\int_{\mR^d}\bar{v}_0(x)\phi(x)\dif x+\int_0^t\int_{\mR^d}\bigg[b_i(s,x)\partial_i\phi(x)+a_{ij}(s,x)\partial_{ij}\phi(x)\no\\
&&+\int_{\mU}\[\phi(x+f(s,x,u))-\phi(x)\]\nu(\dif u)\bigg]v_s(x)\dif x\dif s.
\label{defiequa}
\ee
Take 
$$
\phi(x)=\lambda_n(x)-2n, \quad \lambda_n(x):=n\lambda\left(\frac{\varrho(x)}{n}\right), \quad n\in\mN,
$$
where $\varrho(x)=(1+|x|^2)^{1/2}$ and $\lambda: \mR_+\mapsto\mR_+$ is 
a  twice continuously differentiable, increasing concave function with 
\ce
\lambda(r)=\left\{\begin{array}{ll} r, \quad 0\leq r\leq 1,\\
2, \quad r\geq 2.
\end{array}
\right.
\de
It is easy to see that $\phi\in C_c^2(\mR^d)$. 
Inserting $\phi(x)$ in (\ref{defiequa}), one can obtain that 
\be
\int_{\mR^d}v_t(x)\lambda_n(x)\dif x&=&\int_{\mR^d}\bar{v}_0(x)\lambda_n(x)\dif x+\int_0^t\int_{\mR^d}\bigg[b_i(s,x)\partial_i\lambda_n(x)+a_{ij}(s,x)\partial_{ij}\lambda_n(x)\no\\
&&+\int_{\mU}\[\lambda_n(x+f(s,x,u))-\lambda_n(x)\]\nu(\dif u)\bigg]v_s(x)\dif x\dif s\no\\
&\leq&\int_{\mR^d}\bar{v}_0(x)\lambda_n(x)\dif x+\int_0^t\int_{\mR^d}\bigg[|b_i(s,x)||\partial_i\lambda_n(x)|+|a_{ij}(s,x)||\partial_{ij}\lambda_n(x)|\no\\
&&+\int_{\mU}|\lambda_n(x+f(s,x,u))-\lambda_n(x)|\nu(\dif u)\bigg]v_s(x)\dif x\dif s\no\\
&\leq&C\int_{\mR^d}\bar{v}_0(x)(1+|x|)\dif x+\int_0^t\int_{\mR^d}\bigg[C\varrho(x)|D\lambda_n(x)|+C\varrho(x)^2\|D^2\lambda_n(x)\|\no\\
&&+\int_{\mU}|\lambda_n(x+f(s,x,u))-\lambda_n(x)|\nu(\dif u)\bigg]v_s(x)\dif x\dif s,
\label{soines1}
\ee
where $\lambda_n(x)\leq C(1+|x|)$ is used in the last inequality and $D\lambda_n(x), D^2\lambda_n(x)$ are denoted as the gradient and Hessian matrix of $\lambda_n(x)$, respectively.

Next, by some elementary computations, it holds that
\be
|D\lambda_n(x)|\leq C, \quad |D\lambda_n(x)|\leq C\frac{\lambda_n(x)}{\varrho(x)},\quad \|D^2\lambda_n(x)\|\leq C\frac{\lambda_n(x)}{\varrho(x)^2}.
\label{soines2}
\ee
 We now consider $\int_{\mU}|\lambda_n(x+f(s,x,u))-\lambda_n(x)|\nu(\dif u)$. When $\varrho(x)>2n$, $\lambda_n(x)=2n, \lambda_n(x+f(s,x,u))\leq 2n$, and
\be
\int_{\mU}|\lambda_n(x+f(s,x,u))-\lambda_n(x)|\nu(\dif u)\leq 4n\nu(\mU)=2\lambda_n(x)\nu(\mU).
\label{soines31}
\ee
When $\varrho(x)\leq2n$, $n\frac{\varrho(x)}{n}\leq \lambda_n(x)$, $\frac{\lambda_n(x)}{\varrho(x)}\geq 1$ and
\be
\int_{\mU}|\lambda_n(x+f(s,x,u))-\lambda_n(x)|\nu(\dif u)\leq C\int_{\mU}|f(s,x,u)|\nu(\dif u)\leq C\varrho(x)\leq C\lambda_n(x).
\label{soines32}
\ee
So, (\ref{soines31}) and (\ref{soines32}) yield that 
\be
\int_{\mU}|\lambda_n(x+f(s,x,u))-\lambda_n(x)|\nu(\dif u)\leq C\lambda_n(x).
\label{soines3}
\ee

Finally, combining (\ref{soines2}) and  (\ref{soines3}) with (\ref{soines1}), we have that 
\ce
\int_{\mR^d}v_t(x)\lambda_n(x)\dif x\leq C\int_{\mR^d}\bar{v}_0(x)(1+|x|)\dif x+C\int_0^t\int_{\mR^d}\lambda_n(x)v_s(x)\dif x \dif s.
\de
Thus, the Gronwall inequality implies that
\ce
\int_{\mR^d}v_t(x)\lambda_n(x)\dif x\leq C,
\de
where $C$ is independent of $n, t$. Note that $\lim\limits_{n\rightarrow\infty}\lambda_n(x)=\varrho(x)$. Therefore, by the Fatou lemma, we have
$$
\sup\limits_{t\in[0,T]}\int_{\mR^d}v_t(x)|x|\dif x\leq\sup\limits_{t\in[0,T]}\int_{\mR^d}v_t(x)\varrho(x)\dif x\leq C.
$$
The proof is complete.
\end{proof}

\subsection{The superposition principle for SDEs with jumps and non-local FPEs}
In this subsection, we state a superposition principle between SDEs with jumps and non-local FPEs.

It is well known that for any $\mu_0\in\cP(\mR^d)$, if $\mP$ is a martingale solution to Eq.(\ref{Eq1}) with the initial law $\mu_0$, by simple computation it holds that $(\mP\circ e_t^{-1})$ is a weak solution of Eq.(\ref{FPE1}) 
starting from $\mu_0$. The natural question is whether the converse result is right. The answer is affirmative. The following theorem describes in detail the relationship between martingale solutions to Eq.(\ref{Eq1}) and weak solutions to Eq.(\ref{FPE1}). 

\bt(\cite[Corollary 1.9]{RXZ})\label{supeprin2}
Suppose that (${\bf H}_{b,\sigma}$) and (${\bf H}_f$) hold. 

(i) For any $\mu_0\in\cP(\mR^d)$, the existence of a martingale solution $\mP$ to Eq.(\ref{Eq1}) with the initial law $\mu_0$ is equivalent to the existence of a weak solution $(\mu_t)_{t\in[0,T]}$ to Eq.(\ref{FPE1}) starting from $\mu_0$. Moreover, $\mP\circ e_t^{-1}=\mu_t$ for any $t\in[0,T]$.

(ii) The uniqueness of the martingale solutions $\mP$ to Eq.(\ref{Eq1}) is equivalent to the uniqueness of the weak solutions $(\mu_t)_{t\in[0,T]}$ to Eq.(\ref{FPE1}).  
\et

\br
Theorem \ref{supeprin2} is usually called as a superposition principle.
\er

\section{Limit theorems of SDEs with jumps}\label{add}

In this section, we take $f(t,x,u)=\gamma g(t,x,u)$ for $\gamma\in\mR$. And then Eq.(\ref{Eq1}) becomes 
\be
\dif X_t=b(t,X_t)\dif t+\sigma(t,X_t)\dif B_t+\gamma\int_{\mU}g(t, X_{t-}, u) N(\dif t, \dif u), \qquad t\in[0,T].
\label{Eq5}
\ee
We consider the following sequence of SDEs with jumps: for any $n\in\mN$,
\be
\dif X_t^n=b^n(t, X_t^n)\dif t+\sigma^n(t, X_t^n)\dif B_t+\gamma^n\int_{\mU}g(t, X^n_{t-}, u) N(\dif t, \dif u), \quad t\in[0,T],
\label{Eq6}
\ee
where $b^n: [0,T]\times\mR^d\mapsto\mR^d$, $\sigma^n: [0,T]\times\mR^d\mapsto \mR^{d\times m}$ are Borel measurable functions and $\{\gamma^n\}$ is a real sequence. We study the relationship between martingale solutions of Eq.(\ref{Eq5}) and that of Eq.(\ref{Eq6}) when $b^n\rightarrow b, a^n\rightarrow a, \gamma^n\rightarrow\gamma$ in some sense, where $a^n:=\frac{1}{2}\sigma^n\sigma^{n*}$. 

The first result in the section is the following theorem. 

\bt\label{limit1}
Suppose that $b^n, b, \sigma^n, \sigma$ satisfy (${\bf H}_{b,\sigma}$) uniformly, $\{\gamma^n\}$ is uniformly bounded, $g$ satisfies (${\bf H}_{f}$), and that Eq.(\ref{FPE2}) has a unique weak solution in $\sL_+$. Let $\mu_0(\dif x)=v_0(x)\dif x\in\cP_1(\mR^d)$, and $\mP^n, \mP$ be martingale solutions of Eq.(\ref{Eq6}) and Eq.(\ref{Eq5}) with the initial law $\mu_0$, respectively. Assume that

(i) $b^n\rightarrow b, a^n\rightarrow a$ in $L^1_{loc}([0,T]\times\mR^d)$, $\gamma^n\rightarrow\gamma$ as $n\rightarrow\infty$;

(ii) $\mP^n\circ e_t^{-1}$ is absolutely continuous with respect to the Lebesgue measure on $\mR^d$ and $v^n_t(x)$ denotes the density, i.e., $v^n_t(x):=\frac{(\mP^n\circ e_t^{-1})(\dif x)}{\dif x}$ for any $t\in[0,T]$, and
$$
\sup\limits_{x\in\mR^d}|v^n_t(x)|\leq C_T.
$$

Then $\mP^n\rightarrow\mP$ in $\cP(D_T)$.
\et
\begin{proof}
{\bf Step 1.} We prove that $\{\mP^n\}_{n\in\mN}$ is tight in $\cP(D_T)$.

By Theorem 4.5 in \cite[Page 356]{jjas}, it is sufficient to check that 

(iii) $\lim\limits_{K\rightarrow\infty}\sup\limits_{n}\mP^n\left(\sup\limits_{t\in[0,T]}|w_t|>K\right)=0$,

(iv) For any stopping time $\tau$, it holds that 
$$
\lim\limits_{\t\rightarrow 0}\sup\limits_{n}\sup\limits_{0\leq\tau<\tau+\t\leq T}\mP^n\left(|w_{\tau+\t}-w_{\tau}|\geq N\right)=0, \quad \forall N>0.
$$

Since for any $n\in\mN$, $\mP^n$ is a martingale solution of Eq.(\ref{Eq6}) with the initial law $\mu_0$, Proposition \ref{martweak} yields that there exists a weak solution $\{(\hat{\Omega}^n,\hat{\sF}^n,\hat{\mP}^n;(\hat{\sF}^n_t)_{t\in[0,T]}), (\hat{B}^n,\hat{N}^n,\hat{X}^n)\}$ of Eq.(\ref{Eq6}) with $\cL_{\hat{X}_t^n}=\mP^n\circ e_t^{-1}$ and $\cL_{\hat{X}_0^n}=\mu_0$. So, by Definition \ref{weaksolu}, it holds that for any $ t\in[0,T]$,
\be
\hat{X}_t^n=\hat{X}_0^n+\int_0^tb^n(s, \hat{X}_s^n)\dif s+\int_0^t\sigma^n(s, \hat{X}_s^n)\dif \hat{B}^n_s+\int_0^t\int_{\mU}\gamma^n g(s, \hat{X}_s^n, u)\hat{N}^n(\dif s, \dif u).
\label{nEq7}
\ee
And then the BDG inequality furthermore gives
\ce
\mE^{\hat{\mP}^n}\left(\sup\limits_{s\in[0,t]}|\hat{X}_s^n|\right)&\leq&\mE^{\hat{\mP}^n}|\hat{X}_0^n|+\mE^{\hat{\mP}^n}\left(\int_0^t|b^n(r, \hat{X}_r^n)|\dif r\right)+\mE^{\hat{\mP}^n}\left(\sup\limits_{s\in[0,t]}\left|\int_0^s\sigma^n(r, \hat{X}_r^n)\dif \hat{B}^n_r\right|\right)\\
&&+\mE^{\hat{\mP}^n}\left(\int_0^t\int_{\mU}|\gamma^n||g(s, \hat{X}_s^n, u)|\hat{N}^n(\dif r, \dif u)\right)\\
&\leq&\mu_0(|\cdot|)+\mE^{\hat{\mP}^n}\left(\int_0^t|b^n(r, \hat{X}_r^n)|\dif r\right)+C\mE^{\hat{\mP}^n}\left(\int_0^t\|\sigma^n(r, \hat{X}_r^n)\|^2\dif r\right)^{1/2}\\
&&+\mE^{\hat{\mP}^n}\left(\int_0^t\int_{\mU}|\gamma^n||g(s, \hat{X}_s^n, u)|\nu(\dif u)\dif r\right)\\
&\leq&\mu_0(|\cdot|)+\mE^{\hat{\mP}^n}\left(\int_0^tC(1+|\hat{X}_r^n|)\dif r\right)+C\mE^{\hat{\mP}^n}\left(\int_0^tC(1+|\hat{X}_r^n|)^2\dif r\right)^{1/2},
\de
and 
\ce
\mE^{\hat{\mP}^n}\left(\sup\limits_{s\in[0,t]}\left(1+|\hat{X}_s^n|\right)\right)\leq\mu_0(1+|\cdot|)+C(t+t^{1/2}) \mE^{\hat{\mP}^n}\left(\sup\limits_{s\in[0,t]}\left(1+|\hat{X}_s^n|\right)\right).
\de
By taking $t_0$ with $C(t_0+t_0^{1/2})<1/2$, we obtain that 
$$
\mE^{\hat{\mP}^n}\left(\sup\limits_{s\in[0,t_0]}|\hat{X}_s^n|\right)\leq 2\mu_0(|\cdot|)+1.
$$
On $[t_0,2t_0], [2t_0, 3t_0], \cdots, [[\frac{T}{t_0}]t_0,T]$, by the same way to the above we deduce and conclude that
\be
\mE^{\hat{\mP}^n}\left(\sup\limits_{t\in[0,T]}|\hat{X}_t^n|\right)\leq 2^{[\frac{T}{t_0}]+1}\mu_0(|\cdot|)+2^{[\frac{T}{t_0}]+1}-1,
\label{maxiesti2}
\ee
where $[\frac{T}{t_0}]$ stands for the largest integer no more than $\frac{T}{t_0}$. Thus, it follows from (\ref{maxiesti2}) that $\{\mP^n\}_{n\in\mN}$ satisfies (iii). 

Next, we still consider Eq.(\ref{nEq7}). For any stopping time $\tau$ and $\t>0$ with $0\leq\tau<\tau+\t\leq T$, it holds that
\ce
\hat{X}_{\tau+\t}^n-\hat{X}_{\tau}^n=\int_{\tau}^{\tau+\t}b^n(s,\hat{X}_s^n)\dif s+\int_{\tau}^{\tau+\t}\sigma^n(s, \hat{X}_s^n)\dif \hat{B}^n_s+\int_{\tau}^{\tau+\t}\int_{\mU}\gamma^n g(s, \hat{X}_s^n, u)\hat{N}^n(\dif s, \dif u).
\de
By the similar deduction to the above, we have that
\ce
\mE^{\hat{\mP}^n}|\hat{X}_{\tau+\t}^n-\hat{X}_{\tau}^n|&\leq& \mE^{\hat{\mP}^n}\int_{\tau}^{\tau+\t}C(1+|\hat{X}_s^n|)\dif s+C\mE^{\hat{\mP}^n}\left(\int_{\tau}^{\tau+\t}C(1+|\hat{X}_s^n|)^2\dif s\right)^{1/2}\\
&\leq&C\mE^{\hat{\mP}^n}\left(\sup\limits_{s\in[0,T]}\left(1+|\hat{X}_s^n|\right)\right)(\t+\t^{1/2})\\
&\leq&C\left(2^{[\frac{T}{t_0}]+1}\mu_0(|\cdot|)+2^{[\frac{T}{t_0}]+1}\right)(\t+\t^{1/2}),
\de 
where the last inequality is based on (\ref{maxiesti2}). Thus, by some elementary computations, it holds that $\{\mP^n\}_{n\in\mN}$ satisfies (iv). And then, $\{\mP^n\}_{n\in\mN}$ is relatively weakly compact. That is, there exists a weak convergence subsequence still denoted as $\{\mP^n\}_{n\in\mN}$. 

{\bf Step 2.} We show that a limit point of $\{\mP^n\}_{n\in\mN}$ is $\mP$.

Assume that a limit point of $\{\mP^n\}_{n\in\mN}$ is $\bar{\mP}$. Since $\mP$ is a martingale solution of Eq.(\ref{Eq5}) with the initial law $\mu_0$, and Eq.(\ref{FPE2}) has a unique weak solution, by Theorem \ref{supeprin2} we only need to prove that $\bar{\mP}$ is a martingale solution of Eq.(\ref{Eq5}) with the initial law $\mu_0$. That is, it is sufficient to check that for $0\leq s<t\leq T$ and a bounded continuous $\bar{\cB}_s$-measurable functional $\chi_s: D_T\mapsto \mR$,
\be
\int_{D_T}\left[\phi(w_t)-\phi(w_s)-\int_s^t(\sA_r\phi+\sB_r\phi)(w_r)dr\right]\chi_s(w)\bar{\mP}(\dif w)=0, \quad \forall \phi\in{C_c^2(\mR^d)}.
\label{esti01}
\ee

Next, note that $\mP^n\circ e_t^{-1}\rightarrow\bar{\mP}\circ e_t^{-1}$ in $\cP(\mR^d)$ and $\mP^n\circ e_0^{-1}=\mu_0=\bar{\mP}\circ e_0^{-1}$. Thus, by (ii), there exists a $\bar{v}=(\bar{v}_t)_{t\in[0,T]}$ with $\bar{v}_t(x)\geq0, \forall x\in\mR^d$ and $\int_{\mR^d}\bar{v}_t(x)\dif x=1$ for any $t\in[0,T]$ such that $\bar{\mP}\circ e_t^{-1}(\dif x)=\bar{v}_t(x)\dif x$ and $v^n_t(\cdot)\rightarrow\bar{v}_t(\cdot)$ in $w^*-L^{\infty}(\mR^d)$, where $w^*-L^{\infty}(\mR^d)$ is the dual space of $C_c(\mR^d)$, and $\bar{v}_0(x)=v_0(x)$.
Besides, by (\ref{maxiesti2}) we obtain that
\be
\sup\limits_{n}\sup\limits_{t\in[0,T]}\int_{\mR^d}|x|v^n_t(x)\dif x\leq C.
\label{unimom}
\ee
By suitable approximation, it holds that 
$$
\sup\limits_{t\in[0,T]}\int_{\mR^d}|x|\bar{v}_t(x)\dif x<\infty.
$$
In the following, set $\nu_{t,x}(\dif z):=\nu(\dif g^{-1}(t,x,\cdot)(z))$, and then 
$$
\sB_t\phi(x)=\int_{\mR^d}\[\phi(x+\gamma z)-\phi(x)\]\nu_{t,x}(\dif z).
$$
Thus, based on Lemma 4.1 in \cite{FX} and Lemma \ref{meaapp} below, we know that for any $\e>0$ and the coefficients $b, a$, there exist $\tilde{b}: [0,T]\times\mR^d\mapsto\mR^d$, $\tilde{a}: [0,T]\times\mR^d\mapsto \mS_+(\mR^d)$, where $\mS_+(\mR^d)$ is the set of nonnegative definite symmetric $d\times d$ real matrices, and a family of measures $\tilde{\nu}_{\cdot,\cdot}$ such that 

(v) $\tilde{b}, \tilde{a}$ are continuous and compactly supported;

(vi) for any $\phi\in{C_c^2(\mR^d)}$, $(t,x)\mapsto\tilde{\sB}_t\phi(x)$ is continuous, where $\tilde{\sB}_t\phi(x):=\int_{\mR^d}\[\phi(x+\gamma z)-\phi(x)\]\tilde{\nu}_{t,x}(\dif z)$, and $\sup\limits_{t\in[0,T], x\in\mR^d}|\tilde{\sB}_t\phi(x)|<\infty$;

(vii)
$$
\int_0^T\int_{\mR^d}\left(|b(t,x)-\tilde{b}(t,x)|+\frac{\|a(t,x)-\tilde{a}(t,x)\|}{1+|x|}+|\sB_t\phi(x)-\tilde{\sB}_t\phi(x)|\right)\bar{v}_t(x)\dif x\dif t<\e.
$$
And then the operators with respect to $\tilde{b}, \tilde{a}, \tilde{\nu}_{\cdot,\cdot}$ are denoted as $\tilde{\sA}_t+\tilde{\sB}_t$. 

Now, we treat (\ref{esti01}). Inserting $\tilde{\sA}_t+\tilde{\sB}_t$, one can estimate (\ref{esti01}) to get
\be
&&\left|\int_{D_T}\left[\phi(w_t)-\phi(w_s)-\int_s^t(\sA_r\phi+\sB_r\phi)(w_r)dr\right]\chi_s(w)\bar{\mP}(\dif w)\right|\no\\
&\leq&\left|\int_{D_T}\left[\phi(w_t)-\phi(w_s)-\int_s^t(\tilde{\sA}_r\phi+\tilde{\sB}_r\phi)(w_r)dr\right]\chi_s(w)\bar{\mP}(\dif w)\right|\no\\
&&+\left|\int_{D_T}\left[\int_s^t\left((\tilde{\sA}_r\phi+\tilde{\sB}_r\phi)(w_r)-(\sA_r\phi+\sB_r\phi)(w_r)\right)dr\right]\chi_s(w)\bar{\mP}(\dif w)\right|\no\\
&=:&I_1+I_2.
\label{esti11}
\ee
To deal with $I_1$, we recall that $\mP^n$ is a martingale solution of Eq.(\ref{Eq6}) with the initial law $\mu_0$, which means that 
$$
\int_{D_T}\left[\phi(w_t)-\phi(w_s)-\int_s^t(\sA^n_r\phi+\sB^n_r\phi)(w_r)dr\right]\chi_s(w)\mP^n(\dif w)=0,
$$
where $\sA^n_r+\sB^n_r$ stands for the generator of Eq.(\ref{Eq6}). So, 
\ce
&&\int_{D_T}\left[\phi(w_t)-\phi(w_s)-\int_s^t(\tilde{\sA}_r\phi+\tilde{\sB}_r\phi)(w_r)dr\right]\chi_s(w)\mP^n(\dif w)\\
&=&\int_{D_T}\left[\int_s^t\left((\sA^n_r\phi+\sB^n_r\phi)(w_r)-(\tilde{\sA}_r\phi+\tilde{\sB}_r\phi)(w_r)\right)dr\right]\chi_s(w)\mP^n(\dif w),
\de
and furthermore
\be
&&\left|\int_{D_T}\left[\phi(w_t)-\phi(w_s)-\int_s^t(\tilde{\sA}_r\phi+\tilde{\sB}_r\phi)(w_r)dr\right]\chi_s(w)\mP^n(\dif w)\right|\no\\
&=&\left|\int_{D_T}\left[\int_s^t\left((\sA^n_r\phi+\sB^n_r\phi)(w_r)-(\tilde{\sA}_r\phi+\tilde{\sB}_r\phi)(w_r)\right)dr\right]\chi_s(w)\mP^n(\dif w)\right|\no\\
&\leq&C\int_{D_T}\int_s^t\left|(\sA^n_r\phi+\sB^n_r\phi)(w_r)-(\tilde{\sA}_r\phi+\tilde{\sB}_r\phi)(w_r)\right|dr\mP^n(\dif w)\no\\
&\leq&C\int_s^t\int_{\mR^d}\bigg[|(b_i^n(r,x)-\tilde{b}_i(r,x))\partial_{i}\phi(x)|+|(a^n_{ij}(r,x)-\tilde{a}_{ij}(r,x))\partial_{ij}\phi(x)|\no\\
&&+|\sB^n_r\phi(x)-\sB_r\phi(x)|+|\sB_r\phi(x)-\tilde{\sB}_r\phi(x)|\bigg]v^n_r(x)\dif x\dif r\no\\
&\leq&C\int_s^t\int_{\mR^d}\bigg[|(b_i^n(r,x)-\tilde{b}_i(r,x))\partial_{i}\phi(x)|+|(a^n_{ij}(r,x)-\tilde{a}_{ij}(r,x))\partial_{ij}\phi(x)|\no\\
&&+\int_{\mR^d}|\phi(x+\gamma^n z)-\phi(x+\gamma z)|\nu_{r,x}(\dif z)+|\sB_r\phi(x)-\tilde{\sB}_r\phi(x)|\bigg]v^n_r(x)\dif x\dif r\no\\
&\leq&C\int_s^t\int_{\mR^d}\bigg[|(b_i^n(r,x)-\tilde{b}_i(r,x))\partial_{i}\phi(x)|+|(a^n_{ij}(r,x)-\tilde{a}_{ij}(r,x))\partial_{ij}\phi(x)|\no\\
&&+\|\phi\|_{C_c^2(\mR^d)}\int_{\mR^d}I_{|x+\gamma z|\leq M}|\gamma^n-\gamma||z|\nu_{r,x}(\dif z)+\int_{\mR^d}I_{|x+\gamma z|> M}|\phi(x+\gamma^n z)|\nu_{r,x}(\dif z)\no\\
&&+|\sB_r\phi(x)-\tilde{\sB}_r\phi(x)|\bigg]v^n_r(x)\dif x\dif r,
\label{00}
\ee
where the fact ${\rm supp}(\phi)\subset B_M$ for $M>0$ is applied in the last inequality. As $n\rightarrow\infty$, based on (i) (v) (vi) and $v^n_r(\cdot)\rightarrow\bar{v}_r(\cdot)$ in $w^*-L^{\infty}(\mR^d)$, (\ref{00}) yields that
\be
I_1&\leq&C\int_s^t\int_{\mR^d}\bigg[|(b_i(r,x)-\tilde{b}_i(r,x))\partial_{i}\phi(x)|+|(a_{ij}(r,x)-\tilde{a}_{ij}(r,x))\partial_{ij}\phi(x)|\no\\
&&\qquad +|\sB_r\phi(x)-\tilde{\sB}_r\phi(x)|\bigg]\bar{v}_r(x)\dif x\dif r\no\\
&\leq&C\int_s^t\int_{\mR^d}\bigg[|b(r,x)-\tilde{b}(r,x)|+\|a(r,x)-\tilde{a}(r,x)\|I_{|x|\leq M}\no\\
&&\qquad+|\sB_r\phi(x)-\tilde{\sB}_r\phi(x)|\bigg]\bar{v}_r(x)\dif x\dif r\no\\
&\leq&C\int_s^t\int_{\mR^d}\bigg[|b(r,x)-\tilde{b}(r,x)|+\frac{\|a(r,x)-\tilde{a}(r,x)\|}{1+|x|}\no\\
&&\qquad+|\sB_r\phi(x)-\tilde{\sB}_r\phi(x)|\bigg]\bar{v}_r(x)\dif x\dif r\no\\
&<&C\e,
\label{esti21}
\ee
where we use $I_{|x|\leq M}\leq(1+M)\frac{1}{1+|x|}$ and (vii) in the third and fourth inequality, respectively.

In the following, we treat $I_2$. By the similar deduction to that in (\ref{esti21}), one can obtain that
\be
I_2&\leq&C\int_{D_T}\int_s^t\left|(\tilde{\sA}_r\phi+\tilde{\sB}_r\phi)(w_r)-(\sA_r\phi+\sB_r\phi)(w_r)\right|dr\bar{\mP}(\dif w)\no\\
&\leq&C\int_s^t\int_{\mR^d}\bigg[|(b_i(r,x)-\tilde{b}_i(r,x))\partial_{i}\phi(x)|+|(a_{ij}(r,x)-\tilde{a}_{ij}(r,x))\partial_{ij}\phi(x)|\no\\
&&\qquad+|\sB_r\phi(x)-\tilde{\sB}_r\phi(x)|\bigg]\bar{v}_r(x)\dif x\dif r\no\\
&<&C\e.
\label{esti31}
\ee
Combining (\ref{esti21}) (\ref{esti31}) with (\ref{esti11}), we get that
\ce
\left|\int_{D_T}\left[\phi(w_t)-\phi(w_s)-\int_s^t(\sA_r\phi+\sB_r\phi)(w_r)dr\right]\chi_s(w)\bar{\mP}(\dif w)\right|<C\e.
\de
Letting $\e\rightarrow 0$, we finally have (\ref{esti01}). The proof is complete.
\end{proof}

\bl\label{meaapp}
For any $\e>0$, there is a family of measures $\tilde{\nu}_{\cdot,\cdot}$ such that for any $\phi\in{C_c^2(\mR^d)}$, 

(i) $(t,x)\mapsto\tilde{\sB}_t\phi(x)$ is continuous and $\sup\limits_{t\in[0,T], x\in\mR^d}|\tilde{\sB}_t\phi(x)|<\infty$;

(ii)
$$
\int_0^T\int_{\mR^d}|\sB_t\phi(x)-\tilde{\sB}_t\phi(x)|\bar{v}_t(x)\dif x\dif t<\e.
$$
\el
\begin{proof}
The method comes from \cite[Lemma 3.8]{RXZ}. By \cite[Lemma 14.50, Page 469]{jj}, it holds that there is a measurable function 
$$
h_{t,x} (\t): [0,T]\times\mR^d\times[0,\infty)\mapsto\mR^d\cup\{\infty\},
$$
such that for $t\in[0,T]$ and $x\in\mR^d$
$$
\nu_{t,x}(A)=\int_0^{\infty}I_A(h_{t,x} (\t))\dif \t, \qquad \forall A\in\sB(\mR^d).
$$
And then for any $\phi\in{C_c^2(\mR^d)}$,
$$
\sB_t\phi(x)=\int_{\mR^d}\[\phi(x+\gamma z)-\phi(x)\]\nu_{t,x}(\dif z)=\int_0^{\infty}\[\phi(x+\gamma h_{t,x} (\t))-\phi(x)\]\dif \t.
$$

Next, by the theory of functional analysis, we know that there exists a sequence of measurable functions $\{h^n_{t,x} (\t), n\in\mN\}$ such that for any $\t\geq 0$ and $n\in\mN$, $(t,x)\mapsto h^n_{t,x} (\t)$ 
is continuous with compact support, $|h^n_{t,x} (\t)|\leq |h_{t,x} (\t)|$ and
\ce
\lim_{n\rightarrow\infty}\int_0^T\int_{\mR^d}\int_0^\infty(|h^n_{t,x} (\t)-h_{t,x} (\t)|^2\land 1)\bar{v}_t(x)\dif \t\dif x\dif t=0.
\de
Thus, for any $\e>0$, there exists a $N\in\mN$ such that
\be
\int_0^T\int_{\mR^d}\int_0^\infty(|h^N_{t,x} (\t)-h_{t,x} (\t)|^2\land 1)\bar{v}_t(x)\dif \t\dif x\dif t<\e.
\label{hnh}
\ee
Now, for $t\in[0,T]$ and $x\in\mR^d$, set
$$
\tilde{\nu}_{t,x}(A):=\int_0^{\infty}I_A(h^N_{t,x} (\t))\dif \t, \qquad \forall A\in\sB(\mR^d),
$$
and then $(t,x)\mapsto\tilde{\nu}_{t,x}(A)$ is continuous. So, it holds that for any $\phi\in{C_c^2(\mR^d)}$,
$$
(t,x)\mapsto\tilde{\sB}_t\phi(x)=\int_{\mR^d}\[\phi(x+\gamma z)-\phi(x)\]\tilde{\nu}_{t,x}(\dif z)
$$
is continuous. Besides, note that
$$
|\tilde{\sB}_t\phi(x)|\leq\int_0^\infty|\phi(x+\gamma h^N_{t,x} (\t))-\phi(x)|\dif \t\leq C\int_0^\infty(|\gamma h^N_{t,x} (\t)|\land 1)\dif \t.
$$
Since $h^N_{t,x} (\t)$ has a compact support in $(t,x)$, we have that 
$$
\sup\limits_{t\in[0,T], x\in\mR^d}|\tilde{\sB}_t\phi(x)|<\infty.
$$
Thus, (i) is proved. 

For (ii), we compute that
\ce
|\sB_t\phi(x)-\tilde{\sB}_t\phi(x)|&\leq&\int_0^{\infty}|\phi(x+\gamma h_{t,x} (\t))-\phi(x+\gamma h^N_{t,x} (\t))|\dif \t\\
&\leq&C\int_0^{\infty}\left(I_{B_l}(h_{t,x} (\t))I_{B_{R+|\gamma|l}}(x)+I_{B^c_{l\vee\frac{|x|-R}{|\gamma|}}}(h_{t,x} (\t))\right)\\
&&\qquad\qquad \times(|h^N_{t,x} (\t)-h_{t,x} (\t)|\land 1)\dif \t\\
&\leq&C\left(\int_0^{\infty}\left(I_{B_l}(h_{t,x} (\t))I_{B_{R+|\gamma|l}}(x)+I_{B^c_{l\vee\frac{|x|-R}{|\gamma|}}}(h_{t,x} (\t))\right)\dif \t\right)^{1/2}\\
&&\qquad\qquad \times\left(\int_0^{\infty}(|h^N_{t,x} (\t)-h_{t,x} (\t)|^2\land 1)\dif \t\right)^{1/2}\\
&\leq&C\left(I_{B_{R+|\gamma|l}}(x)\nu_{t,x}(B_l)+\nu_{t,x}(B^c_{l\vee\frac{|x|-R}{|\gamma|}})\right)^{1/2}\\
&&\qquad\qquad \times\left(\int_0^{\infty}(|h^N_{t,x} (\t)-h_{t,x} (\t)|^2\land 1)\dif \t\right)^{1/2}\\
&\leq&C\left(\int_0^{\infty}(|h^N_{t,x} (\t)-h_{t,x} (\t)|^2\land 1)\dif \t\right)^{1/2},
\de
where $l>0$ is a constant, and ${\rm supp}(\phi)\subset B_R$ and $\nu_{t,x}(\mR^d)<\infty$ are used in the second and fifth inequality, respectively. Therefore, the H\"older inequality implies that
\ce
&&\int_0^T\int_{\mR^d}|\sB_t\phi(x)-\tilde{\sB}_t\phi(x)|\bar{v}_t(x)\dif x\dif t\\
&\leq&C\int_0^T\int_{\mR^d}\left(\int_0^{\infty}(|h^N_{t,x} (\t)-h_{t,x} (\t)|^2\land 1)\dif \t\right)^{1/2}\bar{v}_t(x)\dif x\dif t\\
&\leq&CT^{1/2}\left(\int_0^T\int_{\mR^d}\int_0^{\infty}(|h^N_{t,x} (\t)-h_{t,x} (\t)|^2\land 1)\dif \t\bar{v}_t(x)\dif x\dif t\right)^{1/2}\\
&\leq&CT^{1/2}\e^{1/2},
\de
where (\ref{hnh}) is used in the last inequality. The proof is complete.
\end{proof}

\br\label{comfi}
If $\gamma^n=\gamma=0$, $b^n, b, \sigma^n, \sigma$ are uniformly bounded, Theorem \ref{limit1} reduces to \cite[Theorem 3.7]{fi}. Therefore, Theorem \ref{limit1} is more general.
\er

Next, we give an example to explain that $b^n, \sigma^n$ usually happen.

\bx\label{bnan}
Assume that $d=m$, $b: [0,T]\times\mR^d\mapsto\mR^d$, $\sigma: [0,T]\times\mR^d\mapsto \mS_+(\mR^d)$ are continuous and satisfy (${\bf H}_{b,\sigma}$), $\gamma\in\mR$, $g$ satisfies (${\bf H}_{f}$), and that Eq.(\ref{FPE2}) has a unique weak solution in $\sL_+$, and $\mP$ is a martingale solution of Eq.(\ref{Eq5}) with the initial law $\mu_0=v_0(x)\dif x\in\cP_1(\mR^d)$. Set
\ce
&&b_i^n(t,x):=\int_{\mR^d}\varphi_n(x-y)b_i(t,y)\dif y, \\
&&a_{ij}^n(t,x):=\int_{\mR^d}\varphi_n(x-y)a_{ij}(t,y)\dif y, \quad i, j=1, 2, \cdots, d,\\
&&\sigma^n(t,x):=\sqrt{2a^n(t,x)}, \quad a^n(t,x)=(a_{ij}^n(t,x)),\\ 
&&\gamma^n:=\frac{n}{n+1}\gamma,
\de
where $\varphi\in C^{\infty}_c(\mR^d)$ is a nonnegative mollifier with the support in $B_1$ and $\int_{\mR^d}\varphi(x)\dif x=1$, $\varphi_n(x)=n^d\varphi(nx)$. And then we also assume that $\mP^n$ is a martingale solution of the corresponding Eq.(\ref{Eq6}) with the initial law $\mu_0$, $\mP^n\circ e_t^{-1}$ is absolutely continuous with respect to the Lebesgue measure on $\mR^d$ and $v^n_t(x)$ denotes the density, i.e., $v^n_t(x):=\frac{(\mP^n\circ e_t^{-1})(\dif x)}{\dif x}$ for any $t\in[0,T]$, and
$$
\sup\limits_{x\in\mR^d}|v^n_t(x)|\leq C_T.
$$
So, $b^n, \sigma^n, \gamma^n, a^n$ satisfy the conditions of Theorem \ref{limit1}. Thus, it holds that $\mP^n\rightarrow\mP$ in $\cP(D_T)$.

For the readers' convenience, we put the verification of $b^n, \sigma^n, a^n$ satisfying the conditions of Theorem \ref{limit1} in the appendix.
\ex

\medskip

Besides, by checking the proof of Theorem \ref{limit1}, we find that the condition ``Eq.(\ref{FPE2}) has a unique weak solution in $\sL_+$" can be replaced by a weak condition ``Eq.(\ref{Eq5}) has a unique martingale solution". That is, we have the following limit theorem.

\bt
Suppose that $b^n, b, \sigma^n, \sigma$ satisfy (${\bf H}_{b,\sigma}$) uniformly, $\{\gamma^n\}$ is uniformly bounded, $g$ satisfies (${\bf H}_{f}$) and Eq.(\ref{Eq5}) has a unique martingale solution. Let $\mu_0(\dif x)=v_0(x)\dif x\in\cP_1(\mR^d)$ and $\mP^n, \mP$ be martingale solutions of Eq.(\ref{Eq6}) and Eq.(\ref{Eq5}) with the initial law $\mu_0$, respectively. Assume that

(i) $b^n\rightarrow b, a^n\rightarrow a$ in $L^1_{loc}([0,T]\times\mR^d)$, $\gamma^n\rightarrow\gamma$ as $n\rightarrow\infty$;

(ii) $\mP^n\circ e_t^{-1}$ is absolutely continuous with respect to the Lebesgue measure on $\mR^d$ and $v^n_t(x)$ denotes the density, i.e., $v^n_t(x):=\frac{(\mP^n\circ e_t^{-1})(\dif x)}{\dif x}$ for any $t\in[0,T]$, and
$$
\sup\limits_{x\in\mR^d}|v^n_t(x)|\leq C_T.
$$

Then $\mP^n\rightarrow\mP$ in $\cP(D_T)$.
\et

\section{Special cases}\label{spe}

In the section, we analyze some special cases for Eq.(\ref{Eq5}), Eq.(\ref{Eq6}) and Eq.(\ref{FPE2}) and give some concrete and verifiable conditions.

\subsection{The case for $\gamma\neq 0, g(t,x,u)=u$}

In the subsection, we take $\mU\in\sB(\mR^d)$ and $g(t,x,u)=u$ and assume  that for any $p\geq 1$
$$
\int_{\mU}|u|^2(1+|u|)^p\nu(\dif u)<\infty.
$$
Thus, Eq.(\ref{Eq5}), Eq.(\ref{Eq6}) and Eq.(\ref{FPE2}) become
\be
\dif X_t=b(t,X_t)\dif t+\sigma(t, X_t)\dif B_t+\gamma\int_{\mU}u N(\dif t, \dif u), \qquad t\in[0,T], \label{Eq10}\\
\dif X_t^n=b^n(t, X_t^n)\dif t+\sigma^n(t, X_t^n)\dif B_t+\gamma^n\int_{\mU}u N(\dif t, \dif u), \quad t\in[0,T], \label{Eq110}\\
\partial_t v_t=-\partial_i(b_iv_t)+\partial_{ij}(a_{ij}v_t)+\int_{\mU}\[v_t(\cdot-\gamma u)-v_t(\cdot)\]\nu(\dif u).
\label{FPE3}
\ee
The following theorem characterizes the relationship between martingale solutions of (\ref{Eq10}) and ones of (\ref{Eq110}).

\bt\label{fpeweaksolu}
Suppose that $b^n, b, \sigma^n, \sigma$ satisfy (${\bf H}_{b,\sigma}$) uniformly, $\{\gamma^n\}$ are uniformly bounded and that for some $q>1$, $|\triangledown b|\in L^1([0,T], L_{loc}^q(\mR^d)), (\partial_i b_i)^-\in L^1([0,T], L^\infty(\mR^d))$, $\|\triangledown \sigma\|^2\in L^1([0,T], L^\infty(\mR^d))$. Let $\mu_0(\dif x)=v_0(x)\dif x\in\cP_1(\mR^d)$ with
$$
\int_{\mR^d}v_0(x)^r(1+|x|^2)^{(r-1)d}\dif x<\infty
$$
for some $r>1$, and $\mP^n, \mP$ be martingale solutions of Eq.(\ref{Eq110}) and Eq.(\ref{Eq10}) with the initial law $\mu_0$, respectively. Assume that

(i) $b^n\rightarrow b, a^n\rightarrow a$ in $L^1_{loc}([0,T]\times\mR^d)$, $\gamma^n\rightarrow \gamma$ as $n\rightarrow\infty$;

(ii) $\mP^n\circ e_t^{-1}$ is absolutely continuous with respect to the Lebesgue measure on $\mR^d$ and $v^n_t(x)$ denotes the density, i.e., $v^n_t(x):=\frac{(\mP^n\circ e_t^{-1})(\dif x)}{\dif x}$ for any $t\in[0,T]$, and
$$
\sup\limits_{x\in\mR^d}|v^n_t(x)|\leq C_T.
$$

Then $\mP^n\rightarrow\mP$ in $\cP(D_T)$.
\et
\begin{proof}
By the proof of Theorem \ref{limit1}, we only need to prove that Eq.(\ref{FPE3}) has a unique weak solution in $\sL_+$. 

First of all, one can take a complete filtered probability space $(\check{\Omega},\check{\sF},\check{\mP};(\check{\sF}_t)_{t\in[0,T]})$, an $(\check{\sF}_t)$-adapted Brownian motion $(\check{B}_t)$ and an $(\check{\sF}_t)$-adapted Poisson random measure $\check{N}(\dif t, \dif u)$ independent of $(\check{B}_t)$ with the intensity measure $\dif t\nu(\dif u)$. We consider the following equation:
\be
\check{X}_t(x)&=&x+\int_0^t\(b(s,\check{X}_s(x))-\int_{\mU}\gamma u\nu(\dif u)\)\dif s+\int_0^t\sigma(s,\check{X}_s(x))\dif \check{B}_s\no\\
&&+\int_0^t\int_{\mU}\gamma u\tilde{\check{N}}(\dif s\dif u), \quad x\in\mR^d, \quad t\in[0,T],
\label{sde1}
\ee
where $\tilde{\check{N}}(\dif s\dif u):=\check{N}(\dif s\dif u)-\nu(\dif u)\dif s$ is the compensated martingale measure of $\check{N}(\dif s\dif u)$. By \cite[Theorem 4.2]{zhang}, one can obtain that there exists a weak solution $\check{X}_t(x)$ of Eq.(\ref{sde1}) satisfying
\be
\sup\limits_{t\in[0,T]}\check{\mE}\left(\int_{\mR^d}|\psi(\check{X}_t(x))|^{r^*}\rho(\dif x)\right)\leq C\|\psi\|^{r^*}_{L_{\rho}^{r^*}}, \quad \forall \psi\in C_c^\infty(\mR^d),
\label{soluesti}
\ee
where $C$ is independent of $\psi$, $\check{\mE}$ denotes the expectation under the probability measure $\check{\mP}$ and 
$$
\frac{1}{r}+\frac{1}{r^*}=1, \qquad \rho(\dif x)=\frac{1}{(1+|x|^2)^d}\dif x, \quad \|\psi\|^{r^*}_{L_{\rho}^{r^*}}:=\int_{\mR^d}|\psi(x)|^{r^*}\rho(\dif x).
$$

Next, we choose a $\check{\sF}_0$-measurable $d$-dimensional random vector $\check{X}_0$ such that $\check{\mP}\circ\check{X}_0^{-1}=\mu_0$. Thus, $\bar{X}_{\cdot}:=\check{X}_{\cdot}(\check{X}_0)$ solves the following equation
\be
\bar{X}_t&=&\check{X}_0+\int_0^t\(b(s,\bar{X}_s)-\int_{\mU}\gamma u\nu(\dif u)\)\dif s+\int_0^t\sigma(s,\bar{X}_s)\dif \check{B}_s\no\\
&&+\int_0^t\int_{\mU}\gamma u\tilde{\check{N}}(\dif s\dif u).\label{axiequ}
\ee
Moreover, it follows from the H\"older inequality and the Jensen inequality that for any $\psi\in C_c^\infty(\mR^d)$,
\ce
\left|\int_{\mR^d}\psi(x)\cL_{\bar{X}_t}(\dif x)\right|&=&\left|\check{\mE}\psi(\bar{X}_t)\right|=\left|\check{\mE}\left[\check{\mE}[\psi(\check{X}_t(x))|x=\check{X}_0]\right]\right|=\left|\int_{\mR^d}\check{\mE}[\psi(\check{X}_t(x))]v_0(x)\dif x\right|\\
&\leq&\left(\int_{\mR^d}|\check{\mE}[\psi(\check{X}_t(x))]|^{r^*}\rho(\dif x)\right)^{1/r^*}\left(\int_{\mR^d}v_0(x)^r(1+|x|^2)^{rd}\rho(\dif x)\right)^{1/r}\\
&\leq&\left(\check{\mE}\int_{\mR^d}|\psi(\check{X}_t(x))|^{r^*}\rho(\dif x)\right)^{1/r^*}\left(\int_{\mR^d}v_0(x)^r(1+|x|^2)^{(r-1)d}\dif x\right)^{1/r}\\
&\leq&C\|\psi\|_{L_{\rho}^{r^*}},
\de
where the last inequality is based on (\ref{soluesti}). By the theory of functional analysis, we know that there exists a $v_t(\cdot)(1+|\cdot|^2)^d\in L_{\rho}^{r}$ such that
\ce
\int_{\mR^d}\psi(x)v_t(x)(1+|x|^2)^d\rho(\dif x)=\int_{\mR^d}\psi(x)v_t(x)\dif x=\int_{\mR^d}\psi(x)\cL_{\bar{X}_t}(\dif x).
\de
So, \cite[4.25 Problem, Page 325]{ks} gives $\cL_{\bar{X}_t}(\dif x)=v_t(x)\dif x$ for any $t\in[0,T]$. By Proposition \ref{martweak} and Theorem \ref{supeprin2}, it holds that $v$ solves the FPE (\ref{FPE3}) in the distribution sense. Moreover, using  $\mu_0(\dif x)=v_0(x)\dif x\in\cP_1(\mR^d)$ and Proposition \ref{soluinit}, we conclude that $v\in\hat{\sL}_+$, where
\ce
\hat{\sL}_+&:=&\bigg\{v=(v_t)_{t\in[0,T]}: v_t\geq 0, \int_{\mR^d}v_t(x)\dif x=1,  ~\mbox{for any}~ t\in[0,T],\\
&&\quad ~\mbox{and}~\sup\limits_{t\in[0,T]}\left(\int_{\mR^d}v_t(x)^{r}(1+|x|^2)^{(r-1)d}\dif x\right)<\infty, \sup\limits_{t\in[0,T]}\left(\int_{\mR^d}|x|v_t(x)\dif x\right)<\infty\bigg\}.
\de

Finally, note that by \cite[Theorem 4.2]{zhang}, Eq.(\ref{axiequ}) has a pathwise unique weak solution with the initial distribution $\mu_0$ at time $0$. So, for any time $s\in[0,T]$ and $\mu_s(\dif x):=v_s(x)\dif x\in\cP_1(\mR^d)$ with
$$
\int_{\mR^d}v_s(x)^r(1+|x|^2)^{(r-1)d}\dif x<\infty,
$$
by the same deduction to that in \cite[Theorem 4.2]{zhang} one can obtain that Eq.(\ref{axiequ}) has a pathwise unique weak solution with the initial distribution $\mu_s$ at time $s$. From this and Theorem 6.2.3 in \cite{sv}, we know that Eq.(\ref{axiequ}) has a unique martingale solution, which with Theorem \ref{supeprin2} yields that $v$ is unique in $\hat{\sL}_+$. The proof is complete.
\end{proof}

Here we remind that $\sigma$ in the above theorem can be degenerate. Let $\sigma=0$ and then Eq.(\ref{Eq10}) goes into
\be
\dif X_t=b(t,X_t)\dif t+\gamma\int_{\mU}u N(\dif t, \dif u), \qquad t\in[0,T].
\label{Eq111}
\ee
We immediately have the following result.

\bc\label{nocom}
Suppose that $b^n, b, \sigma^n$ satisfy (${\bf H}_{b,\sigma}$) uniformly, $\{\gamma^n\}$ are uniformly bounded and that for some $q>1$, $|\triangledown b|\in L^1([0,T], L_{loc}^q(\mR^d)), (\partial_i b_i)^-\in L^1([0,T], L^\infty(\mR^d))$. Let $\mu_0(\dif x)=v_0(x)\dif x\in\cP_1(\mR^d)$ with
$$
\int_{\mR^d}v_0(x)^r(1+|x|^2)^{(r-1)d}\dif x<\infty
$$
for some $r>1$, and $\mP^n, \mP$ be martingale solutions of Eq.(\ref{Eq110}) and Eq.(\ref{Eq111}) with the initial law $\mu_0$, respectively. Assume that

(i) $b^n\rightarrow b, a^n\rightarrow 0$ in $L^1_{loc}([0,T]\times\mR^d)$, $\gamma^n\rightarrow \gamma$ as $n\rightarrow\infty$;

(ii) $\mP^n\circ e_t^{-1}$ is absolutely continuous with respect to the Lebesgue measure on $\mR^d$ and $v^n_t(x)$ denotes the density, i.e., $v^n_t(x):=\frac{(\mP^n\circ e_t^{-1})(\dif x)}{\dif x}$ for any $t\in[0,T]$, and
$$
\sup\limits_{x\in\mR^d}|v^n_t(x)|\leq C_T.
$$

Then $\mP^n\rightarrow\mP$ in $\cP(D_T)$.
\ec

\br\label{wijupuju}
This corollary means that SDEs with jumps can converge to SDEs with pure jumps in some sense. 
\er

\subsection{The case for $\sigma\neq 0, \gamma=0$}\label{gam}

In the subsection, we take $\sigma\neq 0, \gamma=0$ and require that $\sigma^n, \sigma$ are independent of the space variable $x$. Thus, 
Eq.(\ref{Eq5}), Eq.(\ref{Eq6}) and Eq.(\ref{FPE2}) go into
\be
&&\dif X_t=b(t,X_t)\dif t+\sigma(t)\dif B_t, \qquad t\in[0,T], \label{Eq8}\\
&&\dif X_t^n=b^n(t, X_t^n)\dif t+\sigma^n(t)\dif B_t+\gamma^n\int_{\mU}g(t,x,u) N(\dif t, \dif u), \quad t\in[0,T], \label{Eq88}\\
&&\partial_t v_t=-\partial_i(b_iv_t)+\partial_{ij}(a_{ij}v_t).
\label{FPE4}
\ee
The following proposition describes the relationship between martingale solutions of Eq.(\ref{Eq8}) and that of Eq.(\ref{Eq88}).
 
\bp\label{nojump}
Suppose that $b^n, b, \sigma^n, \sigma, \{\gamma^n\}$ are uniformly bounded, $g$ satisfies (${\bf H}_{f}$), and that $b\in L^1([0,T], BV_{loc}(\mR^d, \mR^d)), \partial_ib_i\in L^1_{loc}([0,T]\times\mR^d),  (\partial_ib_i)^-\in L^1([0,T], L^{\infty}(\mR^d))$. Let $\mu_0(\dif x)=v_0(x)\dif x\in\cP_1(\mR^d)$ with $\|v_0\|_{\infty}<\infty$, and $\mP^n, \mP$ be martingale solutions of Eq.(\ref{Eq88}) and Eq.(\ref{Eq8}) with the initial law $\mu_0$, respectively. Assume that

(i) $b^n\rightarrow b$ in $L^1_{loc}([0,T]\times\mR^d)$, $a^n\rightarrow a$ in $L^1_{loc}([0,T])$, $\gamma^n\rightarrow0$ as $n\rightarrow\infty$;

(ii) $\mP^n\circ e_t^{-1}$ is absolutely continuous with respect to the Lebesgue measure on $\mR^d$ and $v^n_t(x)$ denotes the density, i.e., $v^n_t(x):=\frac{(\mP^n\circ e_t^{-1})(\dif x)}{\dif x}$ for any $t\in[0,T]$, and
$$
\sup\limits_{x\in\mR^d}|v^n_t(x)|\leq C_T.
$$

Then $\mP^n\rightarrow\mP$ in $\cP(D_T)$.
\ep
\begin{proof}
By \cite[Theorem 4.12]{fi} and Proposition \ref{soluinit}, it holds that Eq.(\ref{FPE4}) has a unique weak solution in
\ce
\tilde{\sL}_+&:=&\bigg\{v=(v_t)_{t\in[0,T]}: v_t\geq 0, \int_{\mR^d}v_t(x)\dif x=1,  ~\mbox{for any}~ t\in[0,T],\\
 &&\quad ~\mbox{and}~\sup\limits_{t\in[0,T]}\|v_t(\cdot)\|_{\infty}<\infty, \sup\limits_{t\in[0,T]}\int_{\mR^d}|x|v_t(x)\dif x<\infty\bigg\}.
\de
Note that $\tilde{\sL}_+\subset\sL_+$. Thus, the remaining proof is similar to that of Theorem \ref{limit1}.
\end{proof}

\br\label{wijunoju}
This proposition means that SDEs with jumps can converge to SDEs without jumps in some sense. 
\er

\subsection{The case for $\sigma=0, \gamma=0$}\label{siggam}

In the subsection, we take $\sigma=0, \gamma=0$. Thus, Eq.(\ref{Eq5}) becomes 
\be
\dif X_t=b(t,X_t)\dif t, \qquad t\in[0,T].
\label{Eq12}
\ee
That is, Eq.(\ref{Eq5}) goes into an ordinary differential equation. So, the following proposition presents the relationship between martingale solutions of Eq.(\ref{Eq6}) and ones of Eq.(\ref{Eq12}).

\bp\label{ODE}
Suppose that $b^n, b, \sigma^n, \{\gamma^n\}$ are uniformly bounded, $g$ satisfies (${\bf H}_{f}$), and that $b\in L^1([0,T], BV_{loc}(\mR^d, \mR^d)), \partial_ib_i\in L^1_{loc}([0,T]\times\mR^d),  (\partial_ib_i)^-\in L^1([0,T], L^{\infty}(\mR^d))$. Let $\mu_0(\dif x)=v_0(x)\dif x\in\cP_1(\mR^d)$ with $\|v_0\|_{\infty}<\infty$, and $\mP^n, \mP$ be martingale solutions of Eq.(\ref{Eq6}) and Eq.(\ref{Eq12}) with the initial law $\mu_0$, respectively. Assume that

(i) $b^n\rightarrow b, a^n\rightarrow 0$ in $L^1_{loc}([0,T]\times\mR^d)$, $\gamma^n\rightarrow 0$ as $n\rightarrow\infty$;

(ii) $\mP^n\circ e_t^{-1}$ is absolutely continuous with respect to the Lebesgue measure on $\mR^d$  and $v^n_t(x)$ denotes the density, i.e., $v^n_t(x):=\frac{(\mP^n\circ e_t^{-1})(\dif x)}{\dif x}$ for any $t\in[0,T]$, and
$$
\sup\limits_{x\in\mR^d}|v^n_t(x)|\leq C_T. 
$$

Then $\mP^n\rightarrow\mP$ in $\cP(D_T)$.
\ep

Since the proof of the above proposition is similar to that of Proposition \ref{nojump} and Theorem \ref{limit1}, we omit it.

\br\label{wijuode}
This proposition means that SDEs with jumps can converge to ordinary differential equations in some sense. Moreover, if one take $\gamma^n=0, b^n=b$, and 
\ce
\sigma^n=\frac{1}{\sqrt{n}}\left(\begin{array}{c}1~ 0~ 0 \cdots 0\\
0~ 1~ 0 \cdots 0\\
\vdots~ \vdots~ \ddots \vdots\\
0~ 0~ 0 \cdots 1
\end{array}
\right)_{d\times m},
\de
the above proposition is just right \cite[Corollary 3.9]{fi}.
\er

\section{Appendix}\label{app}

In this section, we show that $b^n, \sigma^n, a^n$ satisfy the conditions in Theorem \ref{limit1}.

{\bf Proof of Remark \ref{bnan}.}

(i) $b^n, \sigma^n$ satisfy (${\bf H}_{b,\sigma}$).

First of all, we show that $b^n$ satisfies (${\bf H}_{b,\sigma}$). For any $t\in[0,T]$, $x\in\mR^d$,
\ce
|b^n(t,x)|&\leq&\int_{\mR^d}\varphi_n(x-y)|b(t,y)|\dif y\leq C\int_{\mR^d}\varphi_n(x-y)(1+|y|)\dif y\\
&\leq&C\int_{\mR^d}\varphi_n(x-y)(1+|x-y|+|x|)\dif y\\
&\leq&C(1+|x|)+C\int_{\mR^d}\varphi_n(x-y)|x-y|\dif y.
\de
Note that 
$$
\int_{\mR^d}\varphi_n(x-y)|x-y|\dif y=\int_{\mR^d}n^d\varphi\(n(x-y)\)|x-y|\dif y=\int_{\mR^d}\varphi(v)\frac{|v|}{n}\dif v\leq1.
$$
Thus, all the above computation yields that
$$
|b^n(t,x)|\leq C(1+|x|),
$$
where $C$ is independent of $n$.

For $\sigma^n$, we only need to verify that $\|a^n(t,x)\|\leq C(1+|x|)^2$. Note that $\|a(t,y)\|\leq C(1+|y|)^2$ and ${\rm supp}(\varphi)=B_1$. The remaining proof is similar to that of $b^n$.

(ii) $b^n\rightarrow b, a^n\rightarrow a$ in $L^1_{loc}([0,T]\times\mR^d)$.

For any $0\leq s<t\leq T$ and $B_R$, 
\ce
&&\int_s^t\int_{B_R}(|b^n(r,x)-b(r,x)|+\|a^n(r,x)-a(r,x)\|)\dif x\dif r\\
&\leq&\int_s^t\int_{B_R}\int_{\mR^d}\varphi_n(x-y)\left(|b(r,y)-b(r,x)|+\|a(r,y)-a(r,x)\|\right)\dif y\dif x\dif r\\
&=&\int_s^t\int_{B_R}\int_{\mR^d}\varphi(v)\left(|b(r,x-\frac{v}{n})-b(r,x)|+\|a(r,x-\frac{v}{n})-a(r,x)\|\right)\dif v\dif x\dif r.
\de
Thus, continuity of $b, \sigma$ in the spatial variable and the dominated convergence theorem allow us to obtain that the above limit is zero as $n\rightarrow \infty$.

\bigskip

\textbf{Acknowledgements:}

The author is very grateful to Professor Renming Song for valuable discussions.

\end{document}